\documentclass{amsart}
\usepackage{amssymb,amsfonts,amsmath}
\usepackage[all]{xy}

\newtheorem{proposition}{Proposition}[section]
\newtheorem{corollary}{Corollary}[section]

\newtheorem{claim}{Claim}[section]
\newtheorem{lemma}{Lemma}[section]
\newtheorem{example}{Example}
\theoremstyle{remark}
\newtheorem{remark}{Remark}
\newcommand{\g}{\gamma}
\newcommand{\si}{\sigma}
\newcommand{\ti}{\tilde}

\newcommand{\n}{\nabla}
\newcommand{\fr}{\frac}

\newcommand{\lf}{\left}
\newcommand{\rg}{\right}
\newcommand{\ep}{\epsilon}
\newcommand{\al}{\alpha}
\newcommand{\real}{\mathbb{R}}

\newcommand{\be}{\beta}
\newcommand{\de}{\delta}

\newcommand{\bl}{\bigl}
\newcommand{\br}{\bigr}

\newcommand{\p}{\partial}

\begin{document}

\title[Splitting, parallel gradient and Bakry-Emery Ricci curvature] {Splitting, parallel gradient  and Bakry-Emery Ricci curvature}
\author{S\'ergio Mendon\c ca}
\address{Departamento de An\'alise, Instituto de Matem\'atica,
Universidade Federal Fluminense, Niter\'oi, RJ, CEP 24020-140,
Brasil} \email{sergiomendonca@id.uff.br}

\dedicatory {To my beloved granddaughters, J\'ulia and Clara}


\begin{abstract} In this paper we obtain a splitting theorem for the symmetric diffusion operator 
$\Delta_\phi=\Delta-\lf<\n\phi,\n \rg>$ and a non-constant $C^3$ function $f$ in a complete Riemannian manifold $M$, under the 
assumptions that the Ricci curvature associated with $\Delta_\phi$ satisfies $\rm Ric_\phi(\n f,\n f)\ge 0$,  that 
$|\n f|$ attains a maximum at $M$ and that $\Delta_\phi$ is non-decreasing along the orbits of $\n f$. 
The proof uses the general fact that a complete manifold $M$ with a non-constant smooth function $f$ with parallel gradient vector field must be a Riemannian product  $M=N\times \mathbb{R}$, where $N$ is any level set of $f$. 
\end{abstract}

\maketitle

\section{\bf Introduction}

Several papers obtained splitting theorems on complete Riemannian manifolds $(M,g)$ 
assuming non-negative sectional curvature, non-negative 
Ricci curvature or non-negative Bakry-Emery Ricci curvature, in the presence of some line in $M$
(see for example [T], [CG], [EH], [FLZ], [WW]). In all these papers the Busemann function $b_\g$ 
associated with a ray $\g$ is 
studied. In general it is proved that the assumptions imply that $b_\g$ is smooth 
and has parallel gradient vector field. In this paper we will not assume the existence of a line.

We will consider the symmetric diffusion operator 
$\Delta_\phi u=\Delta u-\lf<\n\phi,\n u\rg>$, where $\Delta$ is the Laplace-Beltrami operator 
and $\phi$ is a given $C^2$ function on $M$.  
The operator $\Delta_\phi$ is used in probability 
theory, potential theory and harmonic analysis on complete and non-compact Riemannian manifolds. Another important motivation is that, when $\Delta_\phi$ is seen as a symmetric operator in $L^2(M,e^{-\phi}dv_g)$, it   is 
unitarily equivalent to the Schr\"odinger operator $\Delta-\fr 14|\n\phi|^2+
\fr 12\Delta_\phi$ in $L^2(M, dv_g)$, where $dv_g$ is the volume element of $(M,g)$ 
(see for example \cite{d}, \cite{w}, \cite{l}).
 
Let $n$ be the dimension of $M$. For $m\in[n,+\infty]$ we follow \cite{l} and define the $m$-dimensional 
Ricci curvature ${\rm Ric}_{mn}$ associated with the operator $\Delta_\phi$ as follows. Set 
${\rm Ric}_{nn}={\rm Ric}$, where $\rm Ric$ is the usual Ricci curvature. If 
$n<m<\infty$ set 
$${\rm Ric}_{mn}(X,X)={\rm Ric}(X,X)+{\rm Hess}(\phi)(X,X)-\fr{|\lf<\n \phi,X\rg>|^2}
{m-n}.
$$
Finally set 
${\rm Ric}_{\infty n}={\rm Ric}_\phi={\rm Ric}+{\rm Hess}(\phi)$. 

Now we can state our first result:
\begin{theoremA} \label{cor} {\it Let $M$ be a complete connected Riemannian manifold. 
Assume that there exist a $C^3$ function $f$ and a $C^2$ function $\phi$ on $M$ satisfying the following conditions:
\begin{enumerate} 
\item \label{ric} $\rm Ric_{\phi}(\n f,\n f)\ge 0$;
\item \label{bound} $|\n f|$ has a positive global maximum;
\item \label{lap} $\Delta_\phi f$ is non-decreasing along the orbits of $\n f$.
\end{enumerate}
Then $f$ is smooth and $M$ is isometric to the Riemannian product $N\times \real$, where $N$ is any level set of $f$. 
Furthermore  it holds that $\phi$ and $f$ are affine functions on each fiber $\{x\}\times\real$.}
\end{theoremA}

\begin{remark} We will see in Section \ref{examples} that each one of conditions 
(\ref{ric}), (\ref{bound}), (\ref{lap}) is essential in Theorem A. \end{remark}

By a similar proof it can be proved a local version for Theorem A:

\begin{theoremB} \label{cor2} {\it Let $M$ be a connected Riemannian manifold. 
Assume that there exist a $C^3$ function $f$ and a $C^2$ function $\phi$ on $M$ satisfying the following conditions:
\begin{enumerate} 
\item \label{ric} $\rm Ric_{\phi}(\n f,\n f)\ge 0$;
\item \label{bound} $|\n f|$ has a positive global maximum on $M$;
\item \label{lap} $\Delta_\phi f$ is non-decreasing along the orbits of $\n f$.
\end{enumerate}
Then $f$ is smooth and for each point $p\in M$ there exist $\ep>0$ and an open neighborhood $V$ of 
$p$ such that $V=N\times (-\ep,\ep)$, where $N$ is some level set of $f|_V$. 
Furthermore  it holds that $\phi$ and $f$ are affine functions on each fiber $\{x\}\times(-\ep,\ep)$.}
\end{theoremB}

By applying Theorem B to some neighborhood of a point $p$ where $|\n f|$ has a local maximum, 
we obtain:

\begin{corollary} \label{cor3} {\it Let $M$ be a connected Riemannian manifold. 
Assume that there exist a $C^3$ function $f$ and a $C^2$ function $\phi$ on $M$ satisfying the following conditions:
\begin{enumerate} 
\item \label{ric} $\rm Ric_{\phi}(\n f,\n f)\ge 0$;
\item \label{bound} $|\n f|$ has a positive local maximum at some point $p\in M$;
\item \label{lap} $\Delta_\phi f$ is non-decreasing along the orbits of $\n f$ in a 
neighborhood of $p$.
\end{enumerate}
Then there exist $\ep>0$ and an open neighborhood $V$ of 
$p$ such that $f|_V$ is smooth and $V=N\times (-\ep,\ep)$, where $N$ is some level set of $f|_V$. 
Furthermore  it holds that $\phi$ and $f$ are affine functions on each fiber $\{x\}\times(-\ep,\ep)$.}
\end{corollary}

\begin{remark} Since $\rm Ric_{\phi}(\n f,\n f)\ge \rm Ric_{mn}(\n f,\n f)$, for all 
$m\in[n,+\infty]$,  
Theorems A, B and Corollary \ref{cor3} also hold if we replace the condition $\rm Ric_{\phi}(\n f,\n f)\ge 0$ by 
the assumption $\rm Ric_{mn}(\n f,\n f)\ge 0$.
\end{remark}

The main fact that will be used in the proof of Theorem A is the following simple general result, 
which does not require curvature conditions or the existence of lines.

\begin{proposition} \label{th} {\it Let $M$ be a complete connected Riemannian manifold. 
Assume that there exists a non-constant smooth function $f:M\to \real$ such 
that $\n f$ is a parallel vector field. Then $M$ is isometric to $N\times \real$, where $N$ is any level set of $f$. Furthermore $f$ must be an affine function 
on each fiber $\{x\}\times\real$. More precisely, if  $N=f^{-1}(\{c\})$ 
and $|\n f|=C$, the obtained isometry 
$\varphi:N\times \real\to M$ maps each fiber $\{x\}\times \real$ onto the image of the orbit of $\n f$ which 
contains $x$, and it holds that $(f\circ\varphi)(x,t)=c+Ct$.}
\end{proposition}

\section{\bf Functions with parallel gradient vector field}

We recall that a smooth vector field $X$ in $M$ is said to be parallel if for any point $p\in M$, any open neighborhood $U$ of $p$, and any smooth vector field $Y$ in $U$, it holds that 
$\bl(\nabla_YX\br)(p)=0$. 

Proposition \ref{th} could be proved by using the de Rham Decomposition Theorem on the universal cover of $M$ with the induced metric. However, we preferred to present a more elementary proof which just uses the following Berger's extension of Rauch's Comparison Lemma (see for example \cite{ce}).

\begin{lemma}[Berger] \label{berger} Consider complete Riemannian manifolds $W, \ti W$ whose dimensions satisfy $\dim(W)\ge \dim(\ti W)$, a smooth positive function $g:[a,b]\to \real$, unit speed geodesics $\g:[a,b]\to W$, 
 $\ti \g:[a,b]\to \ti W$, and unit parallel vector fields $E$ along $\g$ and $\ti E$ along $\ti\g$, satisfying $\lf<E,\g'\rg>=\bl<\ti E,\ti\g'\br>=0$. 
 For $\de>0$ and $(s,u)\in [a,b]\times[0,\de]$, set 
$\psi_s(u)=\psi^u(s)=\exp_{\g(s)}u g(s) E(s)$ and $\ti\psi_s(u)=\ti\psi^u(s)=\exp_{\ti\g(s)}ug(s)\ti E(s)$. Assume that, for any $s\in [a,b]$, the geodesic 
$\psi_s:[0,\de]\to W$ is free of focal points with respect to 
$\psi_s(0)=\g(s)$. Assume further that, for any $(s,u)\in [a,b]\times[0,\de]$,  
any unit vector $v\in T_{\psi_s(u)}W$ with $\lf<v,\psi_s'(u)\rg>=0$, and 
any unit vector $\ti v\in T_{\ti\psi_s(u)}\ti W$ with $\bl<\ti v,\ti\psi_s'(u)\br>=0$, the sectional curvatures satisfy
$$K\bl(v, \psi_s'(u)\br)=\bl<R\lf(v,\psi_s'(u)\rg)\psi_s'(u),v\br>\ge \bl<\ti R\bl(\ti v,\ti \psi_s'(u)\br)\ti \psi_s'(u),
\ti v\br>=\ti K\bl(\ti v, \ti\psi_s'(u)\br),$$
where $R,\ti R$ are the corresponding tensor curvatures of $W$, respectively, $\ti W$. Then 
it holds that the length $L(\psi^u)\le L\bl(\ti\psi^u\br)$, for any $u\in[0,\de]$.
\end{lemma}

\begin{remark} In the statement of the above Berger's Lemma in \cite{ce}, it was assumed 
that $K(\mu,\nu)\ge \ti K(\ti \mu,\ti \nu)$ for any orthonormal vectors $\mu,\nu\in T_xW$, 
any orthonormal vectors $\ti \mu,\ti \nu\in T_{\ti x}\ti W$ and any $x\in W,\ \ti x\in \ti W$.   
However, the same proof as in \cite{ce} may be used to prove the more 
general formulation as in  Lemma \ref{berger} above. 
\end{remark}

Consider a $C^1$ function $g$ on a manifold such that $|\n g|$ is a constant $D$. Let $\mu$ be an orbit of $\n g$. We recall 
the following simple well-known equality:
\begin{eqnarray}\label{image_real}
g\bl(\mu(t)\br)&=&g\bl(\mu(a)\br)+\int_a^t(g\circ\mu)'(s)ds=g\bl(\mu(a)\br)+\int_a^t\lf<\n g\bl(\mu(s)\br),\mu'(s)\rg>ds\nonumber\\ 
&=&g\bl(\mu(a)\br)+\int_a^t\lf|\n g\bl(\mu(s)\br)\rg|^2 ds=g\bl(\mu(a)\br)+D^2\,(t-a).
\end{eqnarray}

\begin{proof}[\bf Proof of Proposition \ref{th}] Let $M$ be a complete Riemannian manifold and assume that there exists a non-constant
smooth function $f$ such that $\n f$ is parallel. In particular $|\n f|$ is a constant $C>0$. 

\begin{claim}\label{curvature} Fix $p\in M$ and a unit vector field $X$ in a 
neighborhood of $p$ which is orthogonal to $\n f(p)$ at $p$. 
Then the sectional curvature 
$$K\lf(X(p),\fr{\n f(p)}C\rg)=0.
$$
\end{claim}
In fact, since $\n f$ is parallel we have that 
$$\lf(\n_{_X}\n_{_{\n f}}\n f-\n_{_{\n f}}\n_{_X}\n f-
\n_{_{[X,\n f]}}\n f\rg)(p)=0,
$$ 
which proves Claim \ref{curvature}.

From now on we fix a level set $N=f^{-1}(\{c\})\subset M$, for some $c\in \real$. 
\begin{claim} \label{tot_geod} $N$ is a totally geodesic embedded hypersurface.
\end{claim}
Indeed, since $\n f$ has no singularities, the local form of submersions imply that $N$ is a smooth embedded hypersurface. Fix $p\in N$ and a geodesic $\si$ in $M$ satisfying $\si(0)=p$ and $\lf<\n f(\si(0)),\si'(0)\rg>=0$. Since $\n f$ and $\si'$ are parallel vector fields along $\si$,  we obtain that 
$(f\circ \si)'(s)=\lf<\n f(\si(s)),\si'(s)\rg>=\lf<\n f(\si(0)),\si'(0)\rg>=0$ for 
all $s$, hence the image of $\si$ is contained in the level set $N$,  which shows 
that $N$ is totally geodesic. Claim \ref{tot_geod} is proved.

Note that the orbits of $\n f$ intersect $N$ orthogonally and do not intersect each other. 
Furthermore they are geodesics, since $\n f$ is parallel along them. 
In particular the normal exponential map $\exp^\perp:TN^\perp\to M$ is injective. 
It is also surjective, since, for each point $p\in M$, Equation (\ref{image_real}) above implies
that the orbit $\nu$ of $\n f$ which contains $p$ 
satisfies $(f\circ\g)(\real)=\real$, hence $\nu$ must intersect (orthogonally) the level set 
$f^{-1}(\{c\})=N$. We conclude that  $\exp^\perp$ is a diffeomorphism. Thus we define the map
$\varphi:N\times\real\to M$ given by 
$$\varphi(x,t)=\exp_x\fr{t\,\n f(x)}C=\exp^\perp\lf(x,\fr{t\,\n f(x)}C\rg)=
\exp^\perp\lf(x,\fr{t\,\n f(x)}{|\n f(x)|}\rg).
$$ 
Since $\exp^\perp$ is a diffeomorphism, it is easy to see that $\varphi$ is also a diffeomorphism.   

Let $P_t$ denote the parallel transport along the unit speed geodesic $\mu(t)=\varphi(x,t)$. By 
using the fact that $\n f$ is parallel along this geodesic, we obtain that
\begin{equation}\label{parallel_transport} 
\fr{\p\varphi}{\p t}(x,t)=\!\mu'(t)\!=\!P_t\bl(\mu'(0)\br)\!=P_t\lf(\fr{\p\varphi}{\p t}(x,0)\rg)\!=\!P_t\lf(\fr{\n f(x)}C\rg)\!=\!\fr{\n f(\varphi(x,t))}C\, .
\end{equation}
In particular $\mu$ is an orbit of the unit vector field $\n \lf(\fr fC\rg)$. Applying (\ref{image_real}) to the function $g=\fr fC$, we obtain that $\lf(\fr fC\rg)\bl(\varphi(x,t)\br)=\lf(\fr fC\rg)\bl(\mu(t)
\br)=\lf(\fr fC\rg)\bl(\mu(0)\br)+t
=\fr cC+t$, hence
\begin{equation}\label{affine}f\bl(\varphi(x,t)\br)=c+Ct, \mbox{ for any }x\in N\mbox{ and 
any }t\in \real.
\end{equation}
Since $\varphi$ is a diffeomorphism, to prove that $M$ is isometric to $N\times\real$, we just need to prove that $d\varphi_{(x,t)}:T_{(x,t)}(N\times\real)\to T_{\varphi(x,t)}M$ is a linear  isometry for any $(x,t)\in N\times\real$.    To do this, 
we will fix $(x,t)\in N\times\real$ and will consider first the curve 
$\al(s)=(x,t+s)$, which satisfies $\al(0)=(x,t)$ and $|\al'(0)|=1$. Then we will show that 
$|(\varphi\circ\al)'(0)|=1=|\al'(0)|$. We will also consider any unit speed geodesic $\be$ 
orthogonal to $\al$ at $(x,t)$. We will show that $\lf<(\varphi\circ\al)'(0),(\varphi\circ\be)'(0)
\rg>=0=\lf<\al'(0),\be'(0)\rg>$ and $|(\varphi\circ\be)'(0)|=1=|\be'(0)|$.  Then we will 
conclude that  $d\varphi_{(x,t)}:T_{(x,t)}(N\times\real)\to T_{\varphi(x,t)}M$ is a linear isometry and $\varphi$ is an isometry. 

By derivating the equality $(\varphi\circ\al)(s)=\varphi(x,t+s)$ and using (\ref{parallel_transport}) 
we obtain  
\begin{equation}\label{vertical_one}(\varphi\circ\al)'(s)=\fr{\p \varphi}{\p s}(x,t+s)=
\fr{\n f\bl(\varphi(x,t+s)\br)}C=\fr{\n f\bl((\varphi\circ\al)(s)\br)}C.\end{equation} 
In particular it holds that
\begin{equation}\label{vertical}|(\varphi\circ\al)'(0)|=\lf|\fr{\n f\bl(\varphi(x,t)\br)}{C}\rg|=1=|\al'(0)|.
\end{equation} 
Fix $\ep>0$.  Consider a unit 
speed geodesic 
$\be:[-\ep,\ep]\to M$ satisfying $\be(0)=(x,t)=\al(0)$,\ \ $|\be'(0)|=1$ and 
$\lf<\al'(0),\be'(0)\rg>=0$.  Since $\be'(0)$ is tangent to the totally geodesic submanifold $N\times \{t\}$, 
we may write $\be(s)=(\eta(s),t)$ where $\eta:[-\ep,\ep]\to N$ is a geodesic in $N$ 
satisfying $\eta(0)=x$ and $|\eta'(0)|=1$.  
Note that $\eta$ is also a geodesic in  
$M$, since $N$ is totally geodesic by Claim \ref{tot_geod}. By using (\ref{affine}) we obtain that
\begin{equation}\label{fiber}f\bl((\varphi\circ\be)(s)\br)=f\bl(\varphi(\eta(s),t)\br)
=c+Ct.
\end{equation}
As a consequence it holds that
\begin{equation}\label{pre_image}(\varphi\circ\be)\bl([-\ep,\ep]\br)\subset f^{-1}\bl(\{c+Ct\}\br).
\end{equation} 
From (\ref{vertical_one}), (\ref{pre_image})  and the equality $\al(0)=\be(0)=(x,t)$, we 
obtain that 
\begin{eqnarray}\label{cross}\lf<(\varphi\circ\be)'(0), (\varphi\circ\al)'(0)\rg>&=&
\lf<(\varphi\circ\be)'(0), \fr{\n f\bl((\varphi\circ\al)(0)\br)}{C}\rg>\nonumber\\ 
&=&\lf<(\varphi\circ\be)'(0), \fr{\n f\bl((\varphi\circ\be)(0)\br)}{C}\rg>=0=\lf<\be'(0),\al'(0)\rg>.
\end{eqnarray}

Now we consider the unit speed geodesic $\be_{_0}(s)=\bl(\eta(s),0\br)$. 
Let $E$ be the unit parallel vector field along 
$\be_{_0}$ which is orthogonal to $N\times\{0\}$  and satisfies $\bl(\eta(s),u\br)=\exp_{\be_{_0}(s)}uE(s)$, for any $u\in \real$.  Set $\psi_s(u)=\psi^u(s)=\exp_{_{\be_{_0}(s)}}uE(s)$.
 In particular we have that
\begin{equation}\label{exp_ts}\psi_s(t)=\psi^t(s)=\exp_{_{\be_{_0}(s)}}tE(s)=(\eta(s),t)=\be(s).
\end{equation} 
Since $\n f$ is parallel and $N$ is totally 
geodesic, the vector field $\ti E(s)=\fr{\n f(\eta(s))}C$ is a unit parallel vector field along $\eta$ which is orthogonal to $N$. Set $\ti\psi_s(u)=\ti\psi^u(s)=\exp_{_{\eta(s)}}u\ti E(s)=\varphi\bl(\eta(s),u\br)$. Thus we obtain from (\ref{parallel_transport}) that 
\begin{equation}\label{derivat}\ti\psi_s'(u)=\fr{\p\varphi}{\p u}\bl(\eta(s),u\br)=\fr{\n f\bl(\varphi(\eta(s),u)\br)}C=\fr{\n f\bl(\ti\psi_s(u)\br)}C.
\end{equation}

 Note also that 
\begin{equation}\label{ti_exp} \ti\psi_s(t)=\ti\psi^t(s)=\varphi\bl(\eta(s),t\br)=(\varphi\circ\be)(s).\end{equation}

To compare curvatures,  
we fix $s\in[-\ep,\ep]$, $u\ge 0$, and unit vectors $v\in T_{\psi_s(u)}(N\times\real)$
 and $\ti v\in T_{\ti\psi_s(u)}M$ satisfying $\lf<v,\psi_s'(u)\rg>=0=\bl<\ti v,\ti\psi_s'(u)\br>$.  By using the Riemannian product $N\times\real$, Claim \ref{curvature} and Equation (\ref{derivat}), we obtain that
 \begin{equation}\label{comp_K}K\bl(v, \psi_s'(u)\br)=0=
 \ K\lf(\ti v, \fr{\n f\bl(\ti \psi_s(u)\br)}C\rg)=K\bl(\ti v, \ti\psi_s'(u)\br).
\end{equation}
Since $\ti\psi_s$ is a geodesic orthogonal to $N$ and $\exp^\perp:TN^\perp\to M$ is a diffeomorphism, we have that $\ti\psi_s$ is free of focal points to $\ti\psi_s(0)$. 
Similarly we have that $\psi_s$ is free of focal points to $\psi_s(0)$. 
From (\ref{exp_ts}) and (\ref{ti_exp}) we have that $\psi^t=\be$ and $\ti\psi^t=\varphi\circ\be$. 
By using (\ref{comp_K}) we may apply Lemma \ref{berger} with $W=N\times\real$, $\ti W=M$, $g=1$, $\g=\be_{_0}$, 
 $\ti\g=\eta$, obtaining that $L(\be)=L(\psi^t)\le L(\ti\psi^t)=L(\varphi\circ\be)$. 
 We apply this lemma again with 
    $W=M$, $\ti W=N\times\real$, $g=1$, $\g=\eta$, 
 $\ti\g=\be_{_0}$, obtaining that $L(\varphi\circ\be)\le L(\be)$. Varying $\ep>0$ we 
 conclude that
 \begin{equation}\label{horizontal} \lf|(\varphi\circ\be)'(0)\rg|=\lf|\be'(0)\rg|.
 \end{equation}

From (\ref{vertical}), (\ref{cross}), (\ref{horizontal}) we obtain that  $d\varphi_{(x,t)}:T_{(x,t)}(N\times\real)\to T_{\varphi(x,t)}M$ is a linear isometry. Thus 
the diffeomorphism $\varphi:N\times \real\to M$ is an isometry. Furthermore we have 
from (\ref{affine}) that $(f\circ\varphi)(x,t)=c+Ct$, 
hence $f\circ\varphi$ is an affine function on the fiber $\{x\}\times \real$.  
Proposition \ref{th} is proved. 
\end{proof}

A similar proof as above proves the following local version for Proposition \ref{th}:

\begin{proposition} \label{th2}{\it Let $M$ be a connected Riemannian manifold. 
Assume that there exists a non-constant smooth function $f$ on an open subset $U$, such 
that $\n f$ is a parallel vector field on $U$. Then, for each point $p\in U$, there 
exist $\ep>0$ and an open neighborhood 
$V\subset U$ of $p$  such that $V$ is isometric to $N\times (-\ep,\ep)$, 
where $N$ is a smooth level set of $f|_V$. Furthermore  $f$ must be an affine function 
on each fiber $\{x\}\times(-\ep,\ep)$}. 
\end{proposition}

\section{\bf Proof of Theorems A and B} 

To prove Theorems A and B we first recall the famous Bochner formula:
\begin{equation}\label{boch}\fr 12\Delta|\n f|^2=\rm Ric(\n f,\n f)+\n f(\Delta f)+\sum_{i=1}^n|\n_{E_i}\n f|^2,
\end{equation}
where $n$ is the dimension of the Riemannian manifold $M$. 
By using (\ref{boch}), a direct calculation leads easily to the generalized Bochner formula  
below (see [L]):
\begin{equation}\label{boch_2}\fr 12\Delta_\phi|\n f|^2=\rm Ric_\phi(\n f,\n f)+\n f(\Delta_\phi f)+\sum_{i=1}^n|\n_{E_i}\n f|^2.
\end{equation}

Assume that the hypotheses of Theorem A hold. Fix $p\in M$ and a local unit vector field $X$ in 
an open normal ball $B$ centered at $p$. Set $X=E_1$ and construct a local orthonormal frame $E_1,\cdots, E_n$ 
in $B$. By the hypotheses of Theorem A, each parcel on the right side of (\ref{boch_2}) 
is nonnegative, hence $\Delta_\phi|\n f|^2\ge 0$. 
Since $|\n f|$ assumes a global maximum we conclude that $|\n f|$ is constant by the 
maximum principle for elliptic linear operators (see Lemma 2.4 in [FLZ]).  The fact that  $\Delta_\phi|\n f|^2=0$ implies that each parcel on the 
right side of (\ref{boch_2}) vanishes. In particular $\sum_{i=1}^n|\n_{{E_i}}\n f|^2=0$, 
hence $\n_{X}\n f=\n_{{E_1}}\n f=0$. Since $p$ and $X$ were chosen arbitrarily, we obtain that $\n f$ is parallel. 
In particular $f$ is smooth. 
By Proposition \ref{th}, there exists an isometry $\varphi:N\times \real\to M$, where $N$ is some fixed level set of $f$,  
and $f$ is an affine function along each fiber $\{x\}\times \real$. 

By Claim \ref{curvature} 
in the proof of Proposition \ref{th}, we have that {\rm Ric}$(\n f,\n f)=0$. Since {\rm Ric\!}
$_\phi(\n f,\n f)=0$ we obtain that {\rm Hess}$(\phi)(\n f,\n f)=0$. Thus the fact that 
$\n_{\n f}\n f=0$ implies easily that $\n f\bl(\n f(\phi)\br)=0$, hence  $\phi$ is an affine function along any orbit $\varphi\bl(\{x\}\times\real\br)$ of $\n f$. 
The proof of Theorem A is complete.

Theorem B is proved by using Proposition \ref{th2} and proceeding similarly as in the proof of Theorem A.

\section{\label{examples}\bf Examples}

In this section we will see that each one 
of conditions (a), (b), (c) in Theorem A is essential, even in the case that $\phi$ is constant. 

\begin{example} Let $M$ be the hyperbolic $n$-dimensional space and $f$ the Busemann 
function associated to some ray $\g$ in $M$.  We know that any orbit $\si$ of $\n f$ is a 
line containing a ray asymptotic to $\g$. It is also known that $|\n f|=1$ and that 
$\Delta f=n-1$ (see [CM]), hence conditions (b) and (c) in Theorem A hold. Thus condition (a) is essential in Theorem A. 
\end{example} 

\begin{example} Consider a smooth curve $\al:\real\to \real^2\times\{0\}\subset\real^3$ such that $\al(t)=\bl(t,g(t),0\br)$ if $t\in (-1,1)$, where $g$ is an even strict convex nonnegative smooth function satisfying $g(0)=0=g'(0)$ and $\lim_{|t|\to 1}g(t)=1$. Assume further that 
$\al(t)=(1,t,0)$ if $t\ge 1$ and $\al(t)=(-1,|t|,0)$ if $t\le -1$. Let 
$M$ be the smooth surface in $\real^3$ obtained rotating the image of $\al$ around the $y$ axis. 
Clearly the Gauss curvature of $M$ is nonnegative. 
Consider the function $F(x,y,z)=y$ and let $f$ be the restriction of $F$ to $M$. Note 
that  $x^2+z^2\le 1$, and $|\n f(p)|=1$ if $p=(x,y,z)\in M$ and $y\ge 1$. 
Since $|\n f|$ is constant outside a compact set, we have that $|\n f|$ 
attains a maximum at some point of $M$.  Thus $M$ satisfies conditions (a) and (b) 
in Theorem A. This shows that condition (c) is essential in this theorem.
\end{example}

\begin{example} Consider the paraboloid $M\subset \real^3$ given by the equation $z=x^2+y^2$ and the function $F:\real^3\to \real$ given by $F(x,y,z)=z^2$. 
Set $f=F|_M$. Since $M$ has positive Gauss curvature, condition (a) in 
Theorem A holds. In $M-\{(0,0,0)\}$  we consider the coordinates
$\varphi(\rho,\theta)=\bl(\rho\cos\theta,\rho\sin\theta, 
\rho^2\br)$ for $\rho>0$.  On these coordinates we have $\varphi_\rho=(\cos\theta,\sin\theta,2\rho)$ 
and $\n f=\fr{4\rho^3}{1+4\rho^2}\,\varphi_\rho$, hence $\lf|\n f\bl(\varphi(\rho,\theta)\br)\rg|=
\fr{4\rho^3}{\sqrt{1+4\rho^2}}$. As a consequence $\n f$ is unbounded, hence condition 
(b) in Theorem A fails. Now we will see that condition (c) holds. Since the orbit of 
$\n f$ at $(0,0,0)$ is trivial, we just need to check that condition (c) holds in $M-\{(0,0,0)\}$. 
A direct computation leads us to
$$\Delta f=\fr{4\rho^2(3+8\rho^2)}{(1+4\rho^2)^2}+\fr{4\rho^2}{1+4\rho^2}\,.
$$
Thus we obtain that 
\begin{equation}\label{deltaf}\varphi_\rho(\Delta f)=\fr{d}{d\rho}\,(\Delta f)=\fr{8\rho(3+4\rho^2)}{(1+4\rho^2)^3}+\fr{8\rho}{(1+4\rho^2)^2}>0.
\end{equation}
Since  $\n f=\fr{4\rho^3}{1+4\rho^2}\,\varphi_\rho$ we have from (\ref{deltaf}) 
that $\n f(\Delta f)>0$ in $M-\{(0,0,0)\}$, hence $\Delta f$ is non-decreasing along the orbits of $\n f$. Thus condition (c) holds, which shows that condition (b) is essential in Theorem A.
\end{example}

\end{document}